\documentclass[11pt,a4paper]{amsart}
\usepackage{amssymb,amsmath,amsthm}
\usepackage{graphicx}
\usepackage{enumerate}
\usepackage{tikz}

\newtheorem{theorem}{Theorem}[]

\newtheorem{lemma}{Lemma}
\newtheorem{conj}{Conjecture}
\newtheorem{proposition}{Proposition}

\newcommand{\F}{\mathcal{F}}
\newcommand{\LL}{\mathcal{L}}
\newcommand{\G}{\mathcal{G}}
\newcommand{\R}{\mathbb{R}}
\newcommand{\h}{\mathcal{H}}

\newcommand{\al}{\alpha}
\newcommand{\be}{\beta}

\begin{document}

\title{Pairwise intersecting convex sets and cylinders in $\R^3$}

\author{Imre B\'ar\'any}

\begin{abstract} We prove that given a finite collection of cylinders in $\R^3$ with the property that any two them intersect, there is a line intersecting an $\al$ fraction of the cylinders where $\al=\frac 1{14}$. This is a special case of an interesting conjecture.
\end{abstract}


\maketitle

\section{Introduction and main result}\label{introd}

\bigskip
The following beautiful conjecture about finite families of convex bodies in $\R^3$ was raised by Mart\'{i}nez, Rold\'an, and Rubin in~\cite{MRR20}.

\begin{conj}\label{con:basic} There is a constant $\al >0$ such that every finite family $\F$ of convex sets in $\R^3$ whose elements are pairwise intersecting contains a subfamily $\F^*$ of size at least $\al |\F|$ that has a line transversal, that is, there is a line intersecting every element of $\F^*$.
\end{conj}

The condition of pairwise intersecting implies that there is a two-dimen\-sional plane intersecting every set in $\F$. Indeed,
projection of the sets in $\F$ to a line gives a finite family of pairwise intersecting intervals in this line. These intervals have a common
point (for instance by Helly's theorem, one-dimensional version, see~\cite{DGK63}), and the plane whose projection is the common point intersects every set in $\F$.
But we need a line and not a plane intersecting a positive fraction of the sets in $\F$.

\medskip
Projecting the sets in $\F$ to a two-dimensional plane gives a family of convex sets in that plane. Assume that a positive fraction of the triples
of the projected family have a point in common. The two-dimensional version of the fractional Helly theorem of Katchalski and Liu~\cite{KaLiu79} implies that
a positive fraction of the projected sets have a point in common. Thus the line whose projection is this point intersects a positive fraction of the sets in $\F$.
So the conjecture holds under the more restrictive condition that a positive fraction of the triples in $\F$ are intersecting.

\medskip
The conjecture also holds under the condition that the sets in $\F$ are {\sl well-rounded}. This means that there is a number $D$ such for every $A \in \F$ there are concentric balls of radii $r$ and $R$ with $r< R <Dr$ such that $A$ contains the smaller ball and the larger ball contains $A$.

\begin{proposition}\label{prop:2balls} If $\F$ is well-rounded with parameter $D$, then there is a set of lines, $  \LL$, of size at most $32D^2$ such that every $A \in \F$ is intersected by some line in $\LL$.
\end{proposition}

The simple proof is given in the last section.

\medskip
But both conditions (``well-rounded" and ``many triples intersecting'') are too restrictive. For instance, if the family $\F$ consists of (finitely many) coplanar lines with no two parallel, then any two of them intersect. But, assuming general position, there is no intersecting triple in $\F$, and $\F$ is not well-rounded. Yet any line in that plane intersects every set in $\F$. This example shows one of the difficulties that emerge when trying to settle Conjecture~\ref{con:basic}.

\medskip
There are examples of 5 convex sets in $\R^3$ with any two of them intersecting but with no transversal line, that is, no line intersecting all 5 of them. The following example is due to David Darrow~\cite{DD}.
Let $\Delta$ be a regular simplex in $\R^3$, its closed inscribed ball be $D$. Its facets are four triangles $F_1,F_2,F_3,F_4$ considered without their edges.
Each pair $F_i,F_j$ shares an edge. Now to both of them add the same point, one which is very close to one endpoint if their common edge but distinct from the endpoint. Then any two sets in the family
$\F=\{D,F_1,F_2,F_3,F_4\}$ of convex sets have a point in common. But they have no line transversal as one can check easily. So a line can meet at most 4/5 of the sets in $\F$. Repeating $\F$ $m$ times is a family (of size $5m$) of pairwise intersecting convex sets in $\R^3$ and no line can meet more than $4m$  sets in the family.

\medskip
In this short note I'm going to prove a simple special case of the conjecture, namely when $\F$ consists of {\sl cylinders}. A cylinder in $\R^3$ is the Minkowksi sum of a convex set $C$ and a line in $\R^3$.

\begin{theorem}\label{th:cyl} Assume $\F$ is a finite family of cylinders in $\R^3$ having the property that any two cylinders in $\F$ intersect. Then there is a line $L$ intersecting at least $\al |\F|$ elements of $\F$ where $\al=\frac 1{14}$.
\end{theorem}

I find Conjecture~\ref{con:basic} and its companion Conjecture~\ref{con:bip} (to be given in the last section) very interesting. These conjectures are the first and apparently hard instances of a series of analogous conjectures, see ~\cite{MRR20}. They represent a fine and subtle version of Helly type theorems, see the old but excellent survey by Danzer, Gr\"unbaum, and Klee~\cite{DGK63}, or the relevant parts of Matou{\v{s}}ek's comprehensive book~\cite{Mat02}.

\section{Proof of Theorem~\ref{th:cyl}}

\bigskip
Assume $|\F|=n$. Every set in $\F$ is of the form $A=K+L$ where $K$ is a convex body and $L$ is a line in $\R^3$. We suppose, rather for convenience than necessity, that the line $L$ contains the origin. This implies that $K \subset A$.

\medskip
We assume further that $K$ is compact using the method which is usual in Helly type problems. Namely, for every pair $A,B \in \F$ where $A=K+L$ fix a point $z(A,B) \in A\cap B$ and define $K^{new}$ as the convex hull of all $z(A,B)$ with $B \in \F$, $B\ne A$. Then $K^{new}$ is compact, $A^{new}=K^{new}+L\subset A$ is a cylinder again, and in the new system $\{A^{new}: A \in \F \}$ any two sets intersect. Moreover, if a line intersects $\al n$ sets in the new system, then the same line intersects the same sets in the old system because $A^{new}\subset A$.

\begin{figure}[h!]
\centering
\includegraphics[scale=0.9]{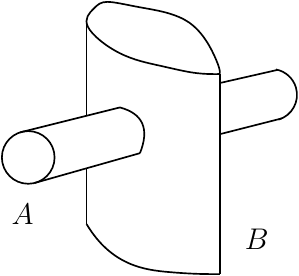}
\caption{A pair representing an $A$ to $B$ arc in $G$.}
\label{fig:AtoB}
\end{figure}

\medskip
We are going to work with an unspecified $\al>0$ that will be fixed later. We define a directed graph $G$ whose vertex set is $\F$,
and $A,B \in \F$ form an arrow $\overrightarrow{AB}$ in $G$ if the lines of $A$ and $B$ are not parallel and $A \setminus B$ is disconnected in $\R^3$,
see Figure~\ref{fig:AtoB}. The outdegree of every $A \in \F$ is at most $\al n$ as otherwise every line in $A$ intersects $\al n$ members of $\F$ and we are done.
Then the average outdegree in $G$ is also at most $\al n$. Consequently at least half of the indegrees
in $G$ is at most $2\al n$.

\medskip
Let $\F'$ be the subfamily of $\F$ consisting of cylinders with indegree at most $2\al n$, so $|\F'|\ge n/2$. Let $C$ be the set in $\F'$ whose width is minimal
and let $P$ be the plane orthogonal to the line of $C$. Define $\F^1$ as the set of cylinders in $\F'$ whose line is parallel with that of $C$.
Observe that the sets $A\cap P$ with $A\in \F^1$ pairwise intersect. So projecting them to any line $L_0$ in $P$ gives pairwise intersecting intervals,
so these intervals have a point in common. Then the line in $P$ containing this point and orthogonal to $L_0$ intersects every set in $\F^1.$ So we are done if $|\F^1|\ge \al n.$

Assume then that $|\F^1|< \al n$. Delete all $A \in \F'$ from $\F'$ with $\overrightarrow{AC}$ an arrow in $G$ together with all $A \in \F^1$ except $C$
to obtain subfamily $\F^* \subset \F'$. Then $C \in \F^*$ and $|\F^*|\ge (\frac 12 -3\al )n$.

\medskip
Project orthogonally every $A \in \F^*$ to $P$. The set $C$ projects to a convex (compact) set $K \subset P$ whose width is $w\ge 0$. Each $A\in \F^*$ projects to a slab $A^*\subset P$
because the line of each $A \in \F^*$ distinct from $C$ is not orthogonal to $P$. The width of $A^*$ is at least $w$ because $C$ has minimal width in $\F^*$.
Note that for each $A \in \F^*$
\begin{enumerate}[{(1)}]
 \item  $A^*\setminus K$ is connected because $\overrightarrow{AC}$ is not an arrow in the digraph $G$,  and
 \item  $A^*\cap K \ne \emptyset$.
\end{enumerate}

\medskip
We consider first the (simple) case $w=0$. Then $K$ is a segment with endpoints $a,b$ in $P$ and $A^* \cap K \ne \emptyset$ because of (2). Further (1)
shows that $A^* \setminus K$ is connected. This together with (2) implies that $A^*$ contains at least one endpoint of the segment $K$ for every $A \in \F^*$.
Consequently one of the lines $a+L$ and $b+L$ intersects half of the sets in $\F^*$, so at least $\frac 12(\frac 12 -3\al)n$ sets from $\F$. So we are done if
$\al \le \frac 12(\frac 12 -3\al)$

\medskip
Thus $w>0$, and we may assume that $w=1$. Write $\h$ for the set of all slabs $S\subset P$ of width at least 1 that intersect but do not cross $K$, that is $S\setminus K$ is connected. In particular, $A^* \in \h$ for every $A \in \F^*$ distinct from $C$. We have now arrived at a statement in plane geometry.

\begin{lemma}\label{l:plane} Under these conditions there is a set consisting of four points that intersects every $S \in \h$.
\end{lemma}

\begin{figure}[h!]
\centering
\includegraphics[scale=0.9]{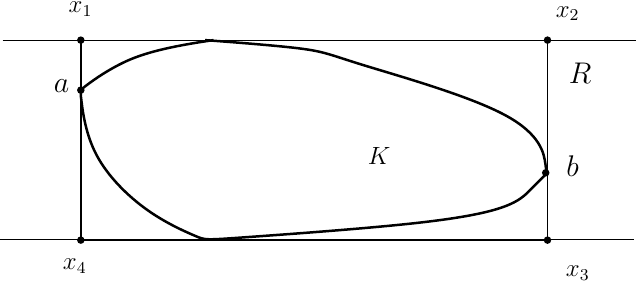}
\caption{Illustration for Lemma~\ref{l:plane}.}
\label{fig:width}
\end{figure}

{\bf Proof.} The minimal slab containing $K$ has width 1, and is bounded by two horizontal (say) lines. Let $x_1,x_2,x_3,x_4$ be the vertices of the
the minimal (with respect to containment) rectangle $R$ that contains $K$ such that the segment $[x_1,x_2]$ is on the upper horizontal line of the slab, and $[x_3,x_4]$ 
is on the lower, see  Figure~\ref{fig:width}.  Let $a$ resp. $b$ be the left-most and right-most (in horizontal direction) points of $K$. Then $a \in [x_1,x_4]$ and $b \in [x_2,x_3]$. 

\medskip
We claim that every $S \in \h$ contains at least one of the points $x_1,x_2,x_3,x_4$. If $S \in \h$ intersects $[x_1,x_4]$, then it contains either $x_1$ 
or $x_4$ because the width of $S$ is at least one, and the segment $[x_1,x_4]$ is of length one. Similarly, if $S$ intersects $[x_2,x_3]$, then it contains either $x_2$ 
or $x_3$. So $S$ is disjoint from the two vertical sides of $R$. Then the endpoints of the segment $[a,b]$ are on different sides of the slab $S$. Then $S$ crosses $K$, so $S \notin \h$.\qed

\medskip
The lemma shows that one of the points $x_i$ $(i=1,2,3,4)$ is contained in the slab $A^*$ for every $A \in \F^*$. So one of the four lines $x_i+L$ intersects $\frac 1{4}(\frac 12-3\al )n$
sets from $\F$. Note that the direction of these lines coincides with the direction with that of $C \in \F$. Requiring $\al=\frac 14\left(\frac 12-3\al \right)$ gives $\al=\frac 1{14}$. This also covers the case $w=0$.

\medskip
{\bf Remarks.} This proof works for the special case of cylinders but does not work in other cases, for instance for half-cylinders that are the Minkowski sum
of a convex sets and a halfline in $\R^3$. The value of $\al$ can be improved to $\frac {7-2\sqrt 6}{25} \approx \frac 1{11.9}$ by choosing $\F'$ to be the subfamily of $\F$ consisting
of cylinders with indegree at most $\beta n$ and optimize the choice of $\beta$. An interesting extension or strengthening
of Theorem~\ref{th:cyl} would be to show that, under the same conditions, there is a set of 100 lines such that each cylinder in $\F$ is intersected by one of them.
Of course, any fixed number instead of 100 would do as well.

\section{The bipartite case}

\bigskip
The paper \cite{MRR20} gives in fact the following, more general conjecture.

\begin{conj}\label{con:bip} There is a constant $\be >0$ such that if  $\F$ and $\G$ are finite families of convex sets in $\R^3$ with the property that $A\cap B\ne \emptyset$ for every pair $A \in \F$ and $B\in \G$, then there is a line intersecting at least either $\be |\F|$ members of $\F$ or $\be |\G|$ elements of $\G$.
\end{conj}

This conjecture implies the one from the introduction. An interesting example is when $\F$ and $\G$ are (finite subsets) of the two sets of lines on the hyperboloid
with equation $x^2+y^2-z^2=1$. The proof method of Theorem~\ref{th:cyl} works in the case when the sets in both families are cylinders but only when their sizes
are (approximately) equal. In this case $\be$ turns out to be equal to the $\al$ from Theorem~\ref{th:cyl}.

\begin{theorem}\label{th:bip} Assume $\F$and $\G$ are finite families of cylinders in $\R^3$ with the property that $A\cap B\ne \emptyset$ whenever $A \in \F$ and $B\in \G$.
If $|\F|=|\G|=n$, then there is a line intersecting at least $\al n$ members of one of $\F$ and $\G$ where $\al=\frac 1{14}$.
\end{theorem}

We only give  a sketch of the {\bf proof} as it is similar to the previous argument. The bipartite and directed graph with classes $\F$ and $\G$, and the arcs
$\overrightarrow{AB}$ and $\overrightarrow{BA}$ are defined the same way as before. The outdegrees are bounded by $\al n$ again otherwise we have a suitable line. The average indegree
in both classes is then at most $\al n$ and we arrive at subfamilies $\F'$ and $\G'$ of size $n/2$ such that the indegrees in each are at most $2\al n$. Again,
$C$ is the minimum width cylinder in $\F' \cup \G'$, say $C \in \F'$. Delete at most $3\al n$ sets from $\G'$, the ones with $\overrightarrow{BC}$ an arc in the directed graph and
the ones whose line is parallel with that of $C$. This way we have a subfamily $\G^*\subset \G'$ whose size is $\left(\frac 12 - 3\al \right)\!n$ and where there is no arc from any $B \in \G^*$ to $C$. An application of Lemma~\ref{l:plane} finishes the proof. \qed

\medskip
{\bf Remark.} The same argument works when $|\F|=n$, $|\G|=m$ and $n \ne m$ but $c< n/m < 1/c$ for some $c>0$. In this case the resulting $\al$ depends on $c$ as well.

\section{Proof of Proposition~\ref{prop:2balls}}

Write $G(a,r)$ for the Euclidean ball of radius $r$ centred at $a \in \R^3$. For every $A \in \F$ we fix a pair of concentric balls $G(a,r)$ and $G(a,R)$ such that  $G(a,r) \subset A \subset G(a,R)$ and $ r<R < Dr$. Let $A_0$ be the convex body in $\F$ for which the radius of the smaller ball is minimal. We can assume that this ball $G(0,1)$, that is, the unit ball centred at the origin. Then $A_0$ is contained in $G(0,D)$.

\medskip
Consider $A \in \F$, $A\ne A_0$ and let $G(a,r)$ and $G(a,R)$ be the pair of concentric balls fixed with $A$. Then $G(0,D)$ intersects $G(a,R)$ since $A_0\cap A \ne \emptyset$.  If $a$ is at distance $T$ from the origin, then $T \le D+R \le D+Dr \le 2rD$ because $r \le 1$ by the choice of $A_0$. Thus $\frac rT > \frac 1 {2D}$. Assume $\ell$ is a line passing through the origin and it makes an angle at most $\varphi = \arcsin \frac 1 {2D}$ with the vector $a$. Such a line must intersect $A$ as it intersects $G(a,r)$.

\medskip
The next (and last) step in the proof is to find a finite set of lines $ \LL$ (whose size is at most $32D^2)$, all passing through the origin so that for every vector $a \ne 0$ there is a line $\ell \in  \LL$ such that the angle between $a$ and $\ell$ is less than $\varphi$. This is a well-known fact, see for instance Lemma 13.1.1 in~\cite{Mat02}. Here is a quick proof: Let $S$ denote the sphere of radius one centred at the origin and for $x \in S$ write $C(x)$ for the set of point $y\in S$ such that the angle between $x$ and $y$ is less than $\varphi/2$. So $C(x)$ is a cap in $S$, and of course,
so is $C(-x)$. Choose points $x_1,x_2,\ldots,x_n$ one by one as long as you can so that the caps $C(x_j)$ and $C(-x_j)$ are disjoint from all previous caps $C(x_i)$ and $C(-x_i)$. The area of such a cap $C(x)$ is at least $\pi \sin^2 (\varphi/2)$ and these $2n$ caps are disjoint, so  $2n \pi \sin^2 (\varphi/2) \le 4\pi$. A simple calculation shows then that $n\le 32D^2$. \qed

\bigskip
{\bf Remark.} This proof works in higher dimensions as well showing that if $\F$ is a finite family of well-rounded (with parameter $D$) convex bodies in $\R^d$, then there is a set of lines, $\LL$, of size $O(D^{d-1})$ such that every $A \in \F$ is intersected by some line in $\LL$. The proof is almost identical with the previous one.

\bigskip
{\bf Acknowledgements.} This piece of work was partially supported by Hungarian National Research Grants No 131529, 131696, and 133819.
My thanks are due to Andreas Holmsen for pointing out an error in a previous version of this paper and suggesting a way to fix it.

\bigskip
\vspace{.5cm} {\sc Imre B\'ar\'any}\\[1mm]
  {\footnotesize R\'enyi Institute of Mathematics}\\[-1mm]
  {\footnotesize 13-15 Re\'altanoda Street, Budapest, 1053 Hungary,}\\[-1mm]
 {\footnotesize e-mail: {\tt barany@renyi.hu}, and }\\[1mm]
  {\footnotesize Department of Mathematics, University College London}\\[-1mm]
  {\footnotesize Gower Street, London WC1E 6BT, UK}\\

\end{document}